\documentclass[12pt]{amsart}
\usepackage{geometry}
\usepackage{graphicx, amsfonts, amssymb}
\usepackage{mathrsfs, amsmath, amsthm}
\usepackage{verbatim}
\usepackage{relsize}
\usepackage{epsfig}

\def\a{\alpha}

\newcommand{\hol}{\rm{Hol}}
\DeclareMathOperator{\og}{O}

\def\D{{\mathbb D}}

\def\C{{\mathbb C}}  
 \def\R{{\mathbb R}}

\def\Dpa{{\mathcal D^p_{\alpha}}}

\def\({\left(}       \def\){\right)}

\DeclareMathOperator{\op}{o} 

\newcommand{\ig}{\stackrel{\text{def}}{=}}

\newtheorem{theorem}{Theorem}

\newtheorem{corollary}{Corollary}
\newtheorem{proposition}{Proposition}
\theoremstyle{definition}

\newtheorem{remark}{Remark}
\numberwithin{equation}{section}
\theoremstyle{theorem}
\newtheorem{other}{\bf Theorem}              
\newtheorem{otherp}{\bf Proposition}  

\newenvironment{pf}{\noindent{\emph{Proof.}}}{$\Box$ }

\begin{document}
\title[Ces\`{a}ro-type operators acting on Dirichlet spaces]
{Ces\`{a}ro-type operators acting on Dirichlet spaces}

\author[O.~Blasco]{\'{O}scar Blasco}
 \address{An\'alisis Matem\'atico,
Universidad de Valencia, 46100, Burjassot, Spain}
\email{Oscar.Blasco@uv.es}
\author[P.~Galanopoulos]{Petros Galanopoulos}
 \address{Department of Mathematics,
Aristotle University of Thessaloniki, 54124, Thessaloniki, Greece}
\email{petrosgala@math.auth.gr}
\author[D.~Girela]{Daniel Girela$^\ast$}
 \address{An\'alisis Matem\'atico,
Universidad de M\'alaga, Campus de Teatinos, 29071 M\'alaga, Spain}
 \email{girela@uma.es}

\begin{abstract} Let $\hol (\mathbb D)$ be the space of all analytic functions in the unit disc $\mathbb D\,=\,\{ z\in \mathbb C : \vert z\vert
<1\}$. For  each $\alpha\in \mathbb R$ we let $\mathcal D^2_\alpha$
be the space of functions  $f\in\hol(\mathbb D)$ such that
$|a_0|^2+\sum_{n=1}^\infty n^{1-\alpha} |a_n|^2<\infty$ where
$f(z)=\sum_{n=0}^\infty a_nz^n$.
\par If $(\eta )=\{ \eta _n\}_{n=0}^\infty $ is a sequence of
complex numbers and $f\in \hol (\mathbb D)$, $f(z)=\sum_{n=0}^\infty
a_nz^n$ ($z\in \mathbb D$), $\mathcal C_{(\eta )}(f)=\mathcal C_{(\{
\eta_n\} )}(f)$ is formally defined by
$$\mathcal C_{(\eta )}(f)=\mathcal C_{\{
\eta_n\}}(f)(z)=\sum_{n=0}^\infty \eta _n\left (\sum
_{k=0}^na_k\right )z^n.$$ The operator $\mathcal C_{(\eta )}$ is a
natural generalization of the Ces\`{a}ro operator. If $\mu $ is a
complex Borel measure on $\mathbb D$ and, for $n=0, 1, 2, \dots $,
$\mu_n=\int_{\mathbb D}w^nd\mu (w)$, the operator $\mathcal C_{\{
\mu _n\} }$ is denoted by $\mathcal C_\mu $.
\par In a recent paper [J. Funct. Anal. \textbf{288} (2025), no. 6, Paper No.~110813], Lin and Xie have studied the question of characterizing the
complex Borel measures $\mu $ on $\mathbb D$ for which the operator
$\mathcal C_\mu$ is bounded (compact) from $\mathcal D^2_\alpha $
into $\mathcal D^2_\beta $ for $\alpha ,\beta >-1$, corresponding to
the spaces of analytic functions $f\in \hol(\mathbb D)$ such that
$f'$ belongs to the Bergman spaces $A^2_\alpha$ and $A^2_\beta$
respectively, and also from $\mathcal D^2_{-1}=S^2$, corresponding
to the space of analytic functions $f\in \hol(\mathbb D)$ such that
$f'\in  H^2$, into itself. They have solved the question for $\alpha
>1$. For the other values of $\alpha $ they have given a number of
conditions which are either necessary or sufficient. They have also
obtained a number of conditions which are either necessary or
sufficient for the boundedness (compactness) of $\mathcal C_\mu $
from $S^2$ into itself.
\par In this paper we give a complete characterization of the
sequences of complex numbers $(\eta_n )$ for which the operator
$\mathcal C_{(\eta )}$ is bounded (compact) from $\mathcal
D^2_\alpha $ into $\mathcal D^2_\beta $ for  $\alpha , \beta \in
\mathbb R$.
\end{abstract}
\thanks{This research has been supported in part by a grant from \lq\lq El Ministerio de
Ciencia e Innovaci\'{o}n\rq\rq , Spain (Project PID2022-138342NB-I00
for the first author and Project PID2022-133619NB-I00  and  grant
from la Junta de Andaluc\'{\i}a FQM-210-G-FEDER for the second and
third authors)}
\thanks{$^\ast$ Corresponding author: Daniel Girela, girela@uma.es} \maketitle


\section{Introduction and preliminaries}\label{intro}
\par\medskip Let $\D=\{z\in\C:|z|<1\}$ denote the open unit disc in the complex plane
$\mathbb C$ and let $\hol(\mathbb D)$ be the space of all analytic
functions in $\D$. Also, $dA$ will denote the area measure on
$\mathbb D$, normalized so that the area of $\mathbb D$ is $1$. Thus
$dA(z)\,=\,\frac{1}{\pi }\,dx\,dy\,=\,\frac{1}{\pi }\,r\,dr\,d\theta
$.\par\medskip The Ces\`{a}ro operator $\mathcal C$ is defined over
the space of all complex sequences as follows: If $(a)=\{ a_k\}
_{k=0}^\infty $ is a sequence of complex numbers then $$\mathcal C
\left ((a)\right )=\left \{ \frac{1}{n+1}\sum_{k=0}^na_k\right \}
_{n=0}^\infty .$$ The operator $\mathcal C$ is known to be bounded
from $\ell^p$ to $\ell^p$ for $1<p\le \infty $. This was proved by
Hardy \cite{H-20} and Landau \cite{L} (see also \cite[Theorem~326,
p.\,\@239~]{HLP}).
\par\medskip
Identifying any given function $f\in \hol (\mathbb D)$ with the
sequence $\{ a_k\}_{k=0}^\infty $ of its Taylor coefficients, the
Ces\`{a}ro operator $\mathcal C$ becomes a linear operator from
$\hol (\mathbb D)$ into itself as follows:
\par If $f\in \hol (\mathbb D)$, $f(z)=\sum_{k=0}^\infty a_kz^k$
\,($z\in \mathbb D$), then
$$\mathcal C(f)(z)=\sum_{n=0}^\infty \left
(\frac{1}{n+1}\sum_{k=0}^na_k\right )z^n,\quad z\in \mathbb D.$$
\par\medskip For $0\,\le\,r\,<\,1$ and
$f$ analytic in $\mathbb D $ we set
$$M_ p(r, f)=\left (
\frac{1}{2\pi }\int_ {0}^{2\pi }\left \vert f(re\sp {i\theta })
\right \vert \sp pd\theta \right )\sp {1/p}, \quad 0<p<\infty , $$
\smallskip
$$M_ \infty (r, f)=\max_ {\vert z\vert =r}\vert g(z)\vert .$$
For $0<p\leq \infty $ the Hardy space $H\sp p$ consists of those
functions $f$, analytic in $\mathbb D $, for which
$$\left \vert \left \vert f \right \vert \right\vert _ {H\sp p}\,\ig\,\sup_ {0<r<1}
M_ p(r, f)<\infty .$$ We refer to \cite{Du:Hp} for the theory of
Hardy spaces.
\par\medskip
For $0<p<\infty $ and $\alpha >-1$ the weighted Bergman space
$A^p_\alpha $ consists of those $f\in \hol (\mathbb D)$ such that
$$\Vert f\Vert _{A^p_\alpha }\,\ig\, \left ((\alpha +1)\int_\mathbb D(1-\vert z\vert ^2)^{\alpha }\vert f
(z)\vert ^p\,dA(z)\right )^{1/p}\,<\,\infty .$$ The unweighted
Bergman space $A^p_0$ is simply denoted by $A^p$. We refer to
\cite{DS,HKZ,Zhu} for the notation and results about Bergman spaces.
\par\medskip
The space $\Dpa$ ($0<p<\infty $, $\alpha >-1$) consists of those
$f\in \hol(\mathbb D)$ such that $f'\in A^p_\alpha $. Hence, if $f$
is analytic in $\mathbb D$, then $f\in \Dpa$ if and only if
\[
\Vert f\Vert _{\Dpa}^p\ig \vert f(0)\vert^p+\Vert f'\Vert
_{A^p_\alpha }^p <\infty .\]
 If $p<\alpha +1$ then it is well known
that $\Dpa =A^p_{\alpha -p}$
 (see, e.~g. Theorem~6 of
\cite{Flett}). On the other hand, if $p>\alpha +2$ then $\Dpa
\subset H^\infty $. Therefore the most interesting ranges for
$(p,\alpha )$ are $\alpha +1\le p\le \alpha +2$.

\par\medskip

The Ces\`{a}ro operator is bounded on $H^p$ for $0<p<\infty $. For
$1<p<\infty $, this follows from a result of Hardy on Fourier series
\cite{H-29} together with the M.~Riesz's theorem on the conjugate
function \cite[Theorem~4.1]{Du:Hp}. Siskakis \cite{Sis-1987} used
semigroups of composition operators to give an alternative proof of
this result and to extend it to $p=1$. A direct proof of the
boundedness on $H^1$ was given by Siskakis in \cite{Sis-1990}. Miao
\cite{Mi} dealt with the case $0<p<1$. The articles \cite{Andersen,
Blasco, Gal, Stem, Sis-1996} include further results on the action
of the Ces\`{a}ro operator on distinct spaces of analytic functions.
\par\medskip
\par\medskip The following generalization of the Ces\`{a}ro operator was introduced in \cite{GGM-2022}. If $\mu $ is a positive finite Borel measure on
$[0,1)$ and $n$ is a non-negative integer, we let $\mu_n$ denote the
moment of order $n$ of $\mu $, that is, $$\mu
_n=\int_{[0,1)}\,t^n\,d\mu (t),\quad n=0, 1, 2, \dots .$$ If $f\in
\hol (\mathbb D)$, $f(z)=\sum_{n=0}^\infty a_nz^n$ ($z\in \mathbb
D$), we define $\mathcal C_\mu (f)$ as follows
$$\mathcal C_\mu (f)(z)\,=\,\sum_{n=0}^\infty \left (\mu_n\sum_{k=0}^n
a_k\right )z^n\,=\,\int_{[0,1)}\frac{f(tz)}{1-tz}\,d\mu (t), \quad
z\in \mathbb D.$$ It is clear that $\mathcal C_\mu $ is a well
defined linear operator $\mathcal C_\mu :\hol (\mathbb D)\rightarrow
\hol (\mathbb D)$. When $\mu $ is the Lebesgue measure on $[0,1)$,
the operator $\mathcal C_\mu $ reduces to the classical Ces\`{a}ro
operator $\mathcal C$.
\par\medskip The work \cite{GGM-2022} includes a complete
characterization of the finite positive Borel measures $\mu $ on
$[0,1)$ for which the operator $\mathcal C_\mu $ is either bounded
or compact on any of the spaces $H^p$ ($1\le p<\infty $),
 $BMOA$, and the Bloch
space $\mathcal B$. A characterization of those $\mu $ for which
$\mathcal C_\mu $ is bounded on $A^p_\alpha $ ($p\ge 1$ $\alpha
>-1$) is also known, it has been obtained in \cite{GGM-2022} for $p>1$ and in \cite{GG-2025}
for $p=1$. The operators $\mathcal C_\mu $ acting on Besov spaces
have been studied in \cite{GGMM}.
\par\medskip Let us turn to deal with two further generalizations of
the Ces\`{a}ro operator. \par (i) If $\mu $ is a complex Borel
measure on $\mathbb D$ (not necessarily supported in $[0,1)$) and
$n\ge 0$, we set
$$\mu_n\,=\,\int_{\mathbb D}w^n\,d\mu (w)$$ and we define the operator $\mathcal C_\mu:\hol (\mathbb D)\rightarrow\hol (\mathbb D)$ as follows:
\par
If $f\in \hol (\mathbb D)$, $f(z)=\sum_{n=0}^\infty a_nz^n$ ($z\in
\mathbb D$), $\mathcal C_\mu (f)$ is defined by
$$\mathcal C_\mu (f)(z)=\sum_{n=0}^\infty\mu_n \left (\sum_{k=0}^na_k\right )z^n=\int_{\mathbb D}\frac{f(wz)}{1-wz}\,d\mu (w),\quad z\in \mathbb D.$$
Since $\{ \mu _n\} $ is a bounded sequence, it is clear that
$\mathcal C_\mu$ is a well defined linear operator from $\hol
(\mathbb D)$ into itself. The complex Borel measures $\mu $ on
$\mathbb D$ for which $C_{\mu }$ is either bounded or compact on
$H^2$ were characterized in \cite{GGM-2023}. This paper contains
also a characterization of those $\mu $ for which $C_{\mu }$ is
bounded on $A^2_\alpha $ ($\alpha >-1$). Lin and Xie \cite{LX} have
characterized those $\mu $ for which $\mathcal C_\mu $ is compact on
$A^2_\alpha $, a question which was left open in \cite{GGM-2023}.
Blasco \cite{Blasco1} has studied the operators $\mathcal C_\mu $
acting acting on the weighted Dirichlet space $\mathcal D_{\rho }$
associated to a sequence of positive numbers $(\rho )=\{\rho_n\} $,
consisting of those $f(z)=\sum_{n=0}^\infty a_nz^n$, analytic in
$\mathbb D$, with $\sum _{n=1}^\infty n\rho_n\vert a_n\vert
^2<\infty $.
\par\medskip (ii) If $(\eta )=\{ \eta _n\}_{n=0}^\infty $ is a sequence of complex
numbers, let $\mathcal C_{(\eta )}=\mathcal C_{\{ \eta_n\}}$ be the
matrix
$$C_{\{ \eta_n\}}\,=\,\left(\begin{array}{ccccccc}
            \eta_0 & 0  & 0  & 0 & . & . \\
  \eta_1 & \eta_1  & 0  & 0 & . & . \\
  \eta_2 & \eta_2  & \eta_2  & 0 & . & . \\
  \eta_3  &  \eta_3 &  \eta_3 &  \eta_3 & . & . \\
   .  & . & . & . & . & . \\.  & . & . & . & . & . \\
\end{array}\right).$$ When $\{ \eta_n\} = \{ \frac{1}{n+1}\}$, $C_{\{
\eta_n\}}$ is the Ces\`{a}ro matrix. Let us remark that this
generalized Ces\`{a}ro matrices are also called Rhaly matrices
\cite{Lei}, \cite{PP} because Rhaly considered particular cases of
them (see \cite{Ra1, Ra2}).\par\medskip The matrix $\mathcal
C_{(\eta )}$ induces an operator, denoted also by $\mathcal C_{(\eta
)}$ or $\mathcal C_{\{ \eta_n\}}$, on the space $\mathcal S$ of
complex sequences by matrix multiplication: If $(a)=\{
a_n\}_{n=0}^\infty $ is a sequence of complex numbers then $$C_{
\{\eta_n\} }(a)\,=\,\left \{ \eta_n\sum_{k=0}^na_n\right
\}_{n=0}^\infty .$$ The sequences $(\eta )=\{ \eta _n\}$ for which
$\mathcal C_{\{ \eta_n\}}$ is bounded on $\ell^p$ ($1<p<\infty $)
have been characterized recently in \cite{GGPr}.
\par Just as above, the operator $\mathcal C_{\{ \eta_n\}}$ can be viewed
as an operator acting on spaces  of analytic functions in the disc.
If $f\in \hol (\mathbb D)$, $f(z)=\sum _{n=0}^\infty a_nz^n $ ($z\in
\mathbb D$), formally, we set
$$\mathcal C_{(\eta )}(f)(z)\,=\,\mathcal C_{\{ \eta _n\}}(f)(z)=\sum_{n=0}^\infty \eta _n\left
(\sum_{k=0}^na_k\right )z^n,\quad z\in\mathbb D$$ whenever the right
hand side makes sense and defines an analytic function in $\mathbb
D$. \par Let us remark that if $\mu $ is a complex Borel measure on
$\mathbb D$ and, for $n\ge 0$,  $\mu_n$ is the moment of order $n$
of $\mu $ then $\mathcal C_{\{ \mu _n\}}\,=\,\mathcal C_\mu .$
\par\medskip The sequence $(\eta )=\{ \eta _n\}$ may not be bounded.
Hence, we cannot be sure that the operator $\mathcal C_{\{
\eta_n\}}$ is a well defined operator from $\hol (\mathbb D)$ into
itself.
\par Notice that for the constant function $1$ we have
$$\mathcal C_{\{
\eta_n\}}(1)\,=\,\sum_{n=0}^\infty \eta_nz^n.$$ Hence, $\mathcal
C_{\{ \eta_n\}}(1)$ is a well defined analytic function in $\mathbb
D$ if and only if the power series $\sum_{n=0}^\infty \eta_nz^n$
defines an analytic function in $\mathbb D$, equivalently, if it has
radius of convergence greater than or equal to $1$.
\par Suppose now that the power series $\sum_{n=0}^\infty \eta_nz^n$
defines an analytic function in $\mathbb D$. This is equivalent to
saying that $\limsup \vert \eta_n\vert ^{1/n}\,\le \,1$. Take $f\in
\hol (\mathbb D)$, $f(z)=\sum_{n=0}^\infty a_nz^n$ ($z\in \mathbb
D$). Set $g(z)=\frac{f(z)}{1-z}$ ($z\in \mathbb D$). The function
$g$ is analytic in $\mathbb D$ and
$$g(z)\,=\,\sum_{n=0}^\infty \left (\sum_{k=0}^na_k\right )z^n,\quad
z\in \mathbb D.$$ Then $\limsup \left \vert
\sum_{k=0}^na_k\right\vert ^{1/n}\,\le \,1$ and it follows that
$$\limsup \left \vert \eta _n\sum_{k=0}^na_n\right \vert ^{1/n}\,\le
1$$ and hence $\mathcal C_{\{ \eta_n\}}(f)$ is a well defined
analytic function in $\mathbb D$.
\par After these considerations we may state the following.
\begin{proposition}\label{prop-welldefined} Let $(\eta )=\{ \eta _n\} _{n=0}^\infty $ be a
sequence of complex numbers.
\begin{itemize}\item[(i)] $\mathcal C_{\{
\eta_n\}}(1)$ is a well defined analytic function in $\mathbb D$ if
and only if the power series $\sum_{n=0}^\infty \eta _nz^n$ defines
an analytic function in $\mathbb D$.
\item[(ii)]  If the power series $\sum_{n=0}^\infty \eta _nz^n$ defines
an analytic function in $\mathbb D$, then $\mathcal C_{\{ \eta_n\}}$
is a well defined linear operator from $\hol (\mathbb D)$ into
itself.
\item[(iii)] If $X$ is a subspace of $\hol (\mathbb D)$ which
contains the constants, then $\mathcal C_{\{ \eta_n\}}$ is a well
defined linear operator from $X$ into $\hol (\mathbb D)$ if and only
if the power series $\sum_{n=0}^\infty \eta _nz^n$ defines an
analytic function in $\mathbb D$.
\end{itemize}
\end{proposition}
\par\medskip From now on, if the power series $\sum_{n=0}^\infty \eta _nz^n$ defines an
analytic function in $\mathbb D$, we shall set
$$F_{(\eta )}(z)\,=\,F_{\{\eta _n\}}(z)\,=\,\sum_{n=0}^\infty \eta
_nz^n,\quad z\in \mathbb D.$$
\par\medskip Lin and Xie have studied in \cite{LX} the question of
characterizing the complex Borel measures $\mu $ on $\mathbb D$ for
which the Ces\`{a}ro-type operator $\mathcal C_\mu $ is bounded or
compact from $\mathcal D^2_\alpha $ into $\mathcal D^2_\beta $ for
distinct values of $\alpha,\beta >-1 $. They have dealt also with
the derivative Hardy space $S^2$ which consists of those $f\in \hol
(\mathbb D)$ such that $f^\prime \in H^2$.
\par The spaces $\mathcal D^2_\alpha $ ($\alpha >-1$) and the space
$S^2$ are Hilbert spaces and it is easy to see if $f\in \hol
(\mathbb D)$, $f(z)=\sum_{n=0}^\infty a_nz^n$ ($z\in \mathbb D$),
then
$$\Vert f\Vert^2 _{\mathcal D^2_\alpha }\,\asymp \,\vert a_0\vert
^2\,+\,\sum_{n=1}^\infty n^{1-\alpha }\vert a_n\vert ^2,$$
$$\Vert f\Vert _{S^2}^2\,=\,\vert a_0\vert ^2\,+\,\sum_{n=1}^\infty
n^2\vert a_n\vert ^2.$$
\par Let us remark that the space $\mathcal D^2_0$ is the Dirichlet
space consisting of those $f\in \hol (\mathbb D)$ with finite
Dirichlet integral.

  We can then unify notation using the equivalent norm and denote by $\mathcal D^2_\alpha$  for any $\alpha\in \R$  the space  of functions $f\in \hol(\mathbb D)$ such that
$$\Vert f\Vert _{\mathcal D^2_\alpha }=( \,\vert a_0\vert
^2\,+\,\sum_{n=1}^\infty n^{1-\alpha }\vert a_n\vert ^2
)^{1/2}<\infty$$ and hence $S^2=\mathcal D^2_{-1}$.

In this paper we shall give a complete description of the sequences
$(\eta_n)$ such that $\mathcal C_{(\eta)} $ is bounded or compact
from $\mathcal D^2_\alpha $ into $\mathcal D^2_\beta $ for distinct
values of $\alpha,\beta \in \mathbb R $. This in particular solves
the questions that were left in \cite{LX} for $\eta_n=\mu_n$ given
by $n$-moments of Borel complex measures $\mu$. Our main results
establish that the situation differs from the cases $\alpha<0$,
$\alpha=0$ and $\alpha>0$  and it can be described looking to the
behaviour of $\|\mathcal C_{(\eta)}(u_N)\|_{\mathcal D^2_\beta}$,
where $u_N(z)=z^N$ for $N\ge 0$. Namely, since $\mathcal
C_{(\eta)}(u_N)(z)=\sum_{n=N}^\infty \eta_nz^n$ we shall prove that
for $\alpha<0$ the boundednees of $\mathcal C_{(\eta)}$ from
$\mathcal D^2_\alpha $ into $\mathcal D^2_\beta $ is actually
equivalent to its compactness and holds if and only if
$F_{(\eta)}\in \mathcal D^2_\beta$ (see Theorem \ref{alphaneg}).
Also we shall see that
 the boundednees (respectively compactness) of $\mathcal C_{(\eta)}$ from $\mathcal D^2_\alpha $ into $\mathcal D^2_\beta $ is
 actually equivalent  to $\|C_{(\eta)}(u_N))\|_{\mathcal D^2_\beta}=O(\frac{1}{\sqrt{\log N}})$ (respectively
$o(\frac{1}{\sqrt{\log N}}$ ) in the case $\alpha=0$ (see Theorem
\ref{alphazero})  and to $\|C_{(\eta)}(u_N))\|_{\mathcal
D^2_\beta}=O(N^{-\alpha/2})$ (respectively $o(N^{-\alpha/2})$ ) in
the case $\alpha>0$ (see Theorem \ref{alphapos}).
\par\medskip

Let us close this section noticing that throughout the paper
 we shall be using the convention that
$C=C(p, \alpha ,q,\beta , \dots )$ will denote a positive constant
which depends only upon the displayed parameters $p, \alpha , q,
\beta \dots $ (which sometimes will be omitted) but not  necessarily
the same at different occurrences. Furthermore, for two real-valued
functions $K_1, K_2$ we write $K_1\lesssim K_2$, or $K_1\gtrsim
K_2$, if there exists a positive constant $C$ independent of the
arguments such that $K_1\leq C K_2$, respectively $K_1\ge C K_2$. If
we have $K_1\lesssim K_2$ and $K_1\gtrsim K_2$ simultaneously, then
we say that $K_1$ and $K_2$ are equivalent and we write $K_1\asymp
K_2$.
\par\medskip
\section{Main results}\label{mainresults}
Let us start by mentioning some previous results which are known on
the above questions in some particular cases:

 Jin and Tang \cite[Theorem 1.2]{JT} give the complete characterization of the boundedness of $C_{(\mu)}$ from $\mathcal D^2_\alpha$ into
 $\mathcal D^2_\beta$ for positive Borel measures $\mu$ defined on $[0,1)$ and for values $0<\alpha,\beta<2$.
 This was extended by Blasco (see \cite[Corollary 6.3]{Blasco1}) for measures $\mu\in M(\mathbb D)$ and all values $\alpha,\beta>0$.

Bao, Guo, Sun, and Wang \cite[Theorem\,\@1.3 and
Theorem\,\@4.3]{BaoGSWa} (see also \cite[Theorem\,\@5]{GG-2025})
have characterized the sequences $(\eta )=\{ \eta _n\} $ for which
the operator $\mathcal C_{(\eta )}$ is bounded from the Dirichlet
space $\mathcal D$ into itself. Theorems 5 and 7 of \cite{GG-2025}
contain also a characterization of the sequences $(\eta )$ for which
$\mathcal C_{(\eta )}$  is either bounded or compact from $\mathcal
D^2_\alpha $ into itself for $0<\alpha \le 1$.
\par Lin and Xie \cite{LX} have obtained a complete characterization
of the complex Borel measures $\mu $ on $\mathbb D$ for which the
operator $\mathcal C_\mu $ is bounded or compact from $\mathcal
D^2_\alpha $ into $\mathcal D^2_\beta $ if $\alpha >1$ and $\beta
>-1$. For the remaining values of $\alpha $, they have obtained a number
of conditions which are either necessary or sufficient for the
boundedness (or compactness) of $\mathcal C_\mu $ from $\mathcal
D^2_\alpha $ into $\mathcal D^2_\beta $. Likewise they give
conditions which are either necessary or sufficient for the
boundedness or compactness of $\mathcal C_\mu $ from $S^2$ into
itself.
\par\medskip
In this paper we obtain a complete characterization of the sequences
$(\eta )=\{ \eta _n\} $ for which the operator $\mathcal C_{(\eta
)}$ is either bounded or compact from $\mathcal D^2_\alpha $ into
$\mathcal D^2_\beta $  for any $\alpha,\beta\in \mathbb R$.
 These
results are summarized in the following theorems. In order to
simplify the notation, if $X$ and $Y$ are two Banach spaces of
analytic functions in $\mathbb D$, $\mathcal B(X, Y)$ will stand for
the space of all bounded linear operators from $X$ into $Y$ and
$\mathcal K(X, Y)$ will stand for the the space of all compact
linear operators from $X$ into $Y$.
\begin{theorem}\label{alphaneg}
Let $(\eta )=\{ \eta _n\}_{n=0}^\infty $ be a sequence of complex
numbers.  If $\alpha <0$ and $\beta \in \mathbb R$, then the
following conditions are equivalent:
\begin{itemize}\item[(a)] $\mathcal C_{(\eta )}\,\in\,\mathcal B(\mathcal
D^2_\alpha , \mathcal D^2_\beta )$.
\item[(b)] $\mathcal C_{(\eta
)}\,\in\,\mathcal K(\mathcal D^2_\alpha , \mathcal D^2_\beta )$.
\item[(c)]
$F_{(\eta )}\,\in\,\mathcal D^2_\beta $.\end{itemize}
\end{theorem}
\begin{corollary}
Let $(\eta )=\{ \eta _n\}_{n=0}^\infty $ be a sequence of complex
numbers.Then
$$\mathcal C_{(\eta )}\,\in\,\mathcal B(\mathcal
S^2 , S^2 ) \Longleftrightarrow \mathcal C_{(\eta )}\,\in\,\mathcal
K(S^2 , S^2 ) \Longleftrightarrow \sum_{n=1}^\infty
n^2|\eta_n|^2<\infty.$$
\end{corollary}

\begin{theorem}\label{alphazero} Let $(\eta )=\{ \eta _n\}_{n=0}^\infty $ be a sequence of complex
numbers. If $\beta \in \mathbb R$, then:
\begin{itemize}\item[(a)] $\mathcal C_{(\eta )}\,\in\,\mathcal B(\mathcal D^2_0 ,
\mathcal D^2_\beta )$\, if and only if
\begin{equation}\label{Obetalog} \sum_{n=N}^\infty \,n^{1-\beta
}\vert \eta_n\vert ^2\,=\,\og \left (\frac{1}{\log N}\right ),\quad
\text{as $N\to \infty $.}\end{equation} \item[(b)] $\mathcal
C_{(\eta )}\,\in\,\mathcal K(\mathcal D^2_0 , \mathcal D^2_\beta
)$\, if and only if
\begin{equation}\label{opbetalog} \sum_{n=N}^\infty \,n^{1-\beta
}\vert \eta_n\vert ^2\,=\,\op \left (\frac{1}{\log N}\right ),\quad
\text{as $N\to \infty $.}\end{equation}\end{itemize}\end{theorem}

\begin{theorem}\label{alphapos}

Let $(\eta )=\{ \eta _n\}_{n=0}^\infty $ be a sequence of complex
numbers. If $\alpha
>0$ and $\beta\in \mathbb R$, then: \begin{itemize} \item[(a)] $\mathcal
C_{(\eta )}\,\in\,\mathcal B(\mathcal D^2_\alpha , \mathcal
D^2_\beta )$\, if and only if
 \begin{equation}\label{Oalphabeta}
\sum_{n=N}^\infty n^{1 -\beta }\vert \eta _n\vert^2\,=\,\og \left
(N^{-\alpha} \right ),\quad \text{as $N\to \infty $.}\end{equation}
 \item[(b)] $\mathcal
C_{(\eta )}\,\in\,\mathcal K(\mathcal D^2_\alpha , \mathcal
D^2_\beta )$\, if and only if
 \begin{equation}\label{oalphabeta}
\sum_{n=N}^\infty n^{1 -\beta }\vert \eta _n\vert^2\,=\,\op \left
(N^{-\alpha} \right ),\quad \text{as $N\to \infty $.}\end{equation}
\end{itemize}
\end{theorem}

\smallskip
\section{Proofs of the main results}\label{proofs}

 We start  recalling the Schur test \cite{Schur} (see also
\cite[p.\,\@9]{Conway} and \cite[Appendix\,\@B]{ARSW}) which will
 be basic in the sequel.\par\medskip
\begin{otherp}[Schur test] Let $\left (\alpha_{j,m}\right )_{j,m=1}^\infty $ be an infinite matrix of complex numbers
and suppose that there exist $c_1, c_2\ge 0$ and a sequence of
positive numbers $\{ p_m\} _{m=1}^\infty $ such that
\begin{align*}\sum_{j=1}^\infty \alpha_{j,m}p_j\,\le
\,c_1p_m,&\quad \text{for all $m$},\\
\sum_{m=1}^\infty \alpha_{j,m}p_m\le c_2p_j,&\quad \text{for all
$j$},\end{align*} then
$$\left \vert \sum_{j,m=1}^\infty \alpha_{j,m}z_jw_m\right \vert ^2\,\le
c_1c_2\left (\sum_{j=1}^\infty \vert z_j\vert ^2\right )\left
(\sum_{j=1}^\infty \vert w_j\vert ^2\right )
$$ for all pair of sequences $\{ z_j\}
_{j=1}^\infty ,\,\{ w_j\} _{j=1}^\infty \,\in \ell^2$.
\end{otherp}

If $X$ is a Banach space of analytic functions in $\mathbb D$, we
shall let $X_0$ be the subspace of $X$ consisting of those $f\in X$
such that $f(0)=0$.
\par Let now $X$ and $Y$ be two Banach spaces of analytic functions in $\mathbb D$ such that $X$ contains the
constants and let $(\eta )=\{ \eta _n\}_{n=0}^\infty $ be a sequence
of complex numbers and assume that $\mathcal C_{(\eta )}\in \mathcal
B(X, Y)$. \par Using Proposition~\ref{prop-welldefined}, we see that
$F_{(\eta )}\in Y$.
\par Also, for any $f\in X$, we have
$$\mathcal C_{(\eta )}(f)\,=\,\mathcal C_{(\eta )}(f(0))\,+\,\mathcal C_{(\eta
)}(f-f(0))\,=\,f(0)F_{(\eta )}\,+\,\mathcal C_{(\eta )}(f-f(0)).$$
Then the following result follows trivially.
\begin{proposition}\label{X0} Let $X$  and $Y$ be two Banach spaces
of analytic functions in $\mathbb D$    such that  $1\in X$  and let
$(\eta )=\{ \eta _n\}_{n=0}^\infty $ be a sequence of complex
numbers. Then $\mathcal C_{(\eta )}\,\in \,\mathcal B(X, Y)$ if and
only if
\begin{equation}\label{X-X0}F_{(\eta )}\,\in Y\quad \text{and \,\,$\mathcal C_{(\eta
)}\in \mathcal B(X_0, Y)$}.\end{equation}
\end{proposition}
\par\bigskip
\centerline {\it Proof of Theorem\,\@\ref{alphaneg}.}\par\medskip
The implication (a)\, $\Rightarrow $\, (c) follows using
Proposition\,\@\ref{X0} and the implication (b)\, $\Rightarrow $\,
(a) is obvious. Hence it remains to prove that (c)\, $\Rightarrow
$\, (b). So assume (c), that is, $\sum_{n=1}^\infty n^{1-\beta
}\vert \eta _n\vert ^2<\infty $. Set
$$A_N\,=\,\sum_{n=N}^\infty n^{1-\beta
}\vert \eta _n\vert ^2, \quad N=1, 2, 3, \dots .$$ Then $\{
A_N\}_{N=1}^\infty $ is a decreasing sequence of nonnegative numbers
and $\{ A_N\} \to 0$, as $N\to \infty $. \par For $N\,=\,1, 2, 3,
\dots $, let $\mathcal C_N:\mathcal D^2_\alpha \rightarrow \mathcal
D^2_\beta $ be the operator defined as follows: If $f\in \mathcal
D^2_{\alpha }$, $f(z)\,=\,\sum_{n=0}^\infty a_nz^n$ ($z\in \mathbb
D$),
$$\mathcal C_N(f)(z)\,=\,\sum_{n=0}^N\eta _n\left
(\sum_{k=0}^na_k\right )z^n,\quad z\in \mathbb D.$$ The operators
$\mathcal C_N$ are finite rank operators from $\mathcal D^2_\alpha $
into $\mathcal D^2_\beta $. Then (b) will follow from the fact that
\begin{equation}\label{fin-rank}\mathcal C_N\,\rightarrow \,\mathcal
C_{(\eta )},\quad\text{in the operator norm.}\end{equation} So, take
$f\in \mathcal D^2_\alpha $, $f(z)=\sum_{n=0}^\infty a_nz^n$ ($z\in
\mathbb D$). We have
\begin{equation}\label{CNC}\Vert \mathcal C_N(f)\,-\,\mathcal
C_{(\eta )}(f)\Vert _{\mathcal D^2_\beta }^2\,= \sum_{n=N+1}^\infty
n^{1-\beta }\vert \eta _n\vert ^2\left (\sum_{k=0}^n\vert a_k\vert
\right )^2\,\lesssim \,I\,+\,II,\end{equation} where
$$I\,=\,\vert a_0\vert ^2A_{N+1},\quad II\,=\,\sum_{n=N+1}^\infty n^{1-\beta }\vert \eta _n\vert ^2\left
(\sum_{k=1}^n\vert a_k\vert \right )^2.$$ Clearly,
$$I\,=\,\vert a_0\vert ^2A_{N+1}\,\le A_N\Vert f\Vert_{\mathcal
D^2_\alpha }^2.$$ Now, using the Cauchy-Schwarz inequality and
bearing in mind that $\alpha <0$, we see that
\begin{align*}II\,&=\,\,\sum_{n=N+1}^\infty n^{1-\beta }\vert \eta
_n\vert ^2\left (\sum_{k=1}^n\vert a_k\vert \right )^2\,=\,
\sum_{n=N+1}^\infty n^{1-\beta }\vert \eta _n\vert ^2\left
(\sum_{k=1}^n \left
(k^{(1-\alpha )/2}\vert a_k\vert \right )k^{(\alpha -1)/2}\right )^2\\
\,&\le\,\,\sum_{n=N+1}^\infty n^{1-\beta }\vert \eta _n\vert ^2\left
(\sum_{k=1}^nk^{1-\alpha }\vert a_k\vert ^2\right
)\left(\sum_{k=1}^\infty k^{\alpha -1}\right)\,\lesssim A_N\Vert
f\Vert _{\mathcal D^2_\alpha }^2.\end{align*} Using these estimates
in (\ref{CNC}), we see that
$$\Vert \mathcal C_N(f)\,-\,\mathcal C_{(\eta )}(f)\Vert _{\mathcal
D^2_\beta }^2\,\lesssim \,A_N\Vert f\Vert _{\mathcal D^2_\alpha
}^2$$ and then (\ref{fin-rank}) follows.

\par\bigskip
\centerline{\it Proof of Theorem\,\@\ref{alphazero} \,\@(a).}\par
\par\medskip
Proposition~\ref{X0} implies that $\mathcal C_{(\eta )}\in \mathcal
B(\mathcal D^2_0, \mathcal D^2_\beta )$ if and only if
$\sum_{n=1}^\infty n^{1-\beta }\vert \eta_n\vert ^2<\infty $ and
$\mathcal C_{(\eta )}\in \mathcal B(\mathcal D^2_{0,0}, \mathcal
D^2_\beta )$. Then it follows that  part (a) of
Theorem\,\@\ref{alphazero} is equivalent to the following.
\begin{theorem}\label{D200beta} Let $(\eta )=\{
\eta_n\}_{n=0}^\infty $ be a sequence of complex numbers and $\beta
\in \mathbb R$. Then $\mathcal C_{(\eta )}\in \mathcal B(\mathcal
D^2_{0,0}, \mathcal D^2_\beta )$ if and only if
\begin{equation}\label{betalog}\sum_{n=N}^\infty n^{1-\beta }\vert \eta_n\vert^2\,=\,\og \left
(\frac{1}{\log N}\right ),\quad \text{as $N\to \infty
$}.\end{equation}
\end{theorem}
\par\medskip
\begin{pf} Assume that $\mathcal C_{(\eta )}\in \mathcal B(\mathcal D^2_{0,0},
\mathcal D^2_\beta )$. For $\frac{1}{2}<b<1$, set
$$h_b(z)\,=\,\left (\log \frac{1}{1-b}\right )^{-1/2}\log
\frac{1}{1-bz}\,=\,\left (\log \frac{1}{1-b}\right
)^{-1/2}\sum_{n=1}^\infty \frac{b^n}{n}z^n,\quad z\in \mathbb D.$$
We have that $h_b\in \mathcal D^2_{0,0}$ and $\Vert h_b\Vert _
{\mathcal D^2_{0,0}}^2\asymp 1$. Then we have, for each $N\ge 2$
\begin{align*} 1\gtrsim\, &\,\Vert \mathcal C_{(\eta )}(h_b)\Vert
_{\mathcal D^2_\beta }^2\,\asymp \,\left (\log \frac{1}{1-b}\right
)^{-1}\sum _{n=1}^\infty n^{1-\beta }\vert \eta _n\vert ^2\left
(\sum_{k=1}^n\frac{b^k}{k}\right )^2\\
\gtrsim\, & \,\left (\log \frac{1}{1-b}\right )^{-1}\left(\sum
_{n=N}^\infty n^{1-\beta }\vert \eta _n\vert ^2\right)\left
(\sum_{k=1}^N\frac{b^k}{k}\right )^2\\ \gtrsim \,&\, \left (\log
\frac{1}{1-b}\right )^{-1}\left(\sum _{n=N}^\infty n^{1-\beta }\vert
\eta _n\vert ^2\right) b^{2N}\log^2 N.
\end{align*}
Taking $b=1-\frac{1}{N}$, we obtain (\ref{betalog}).
\par\medskip Let us prove the other implication. So, assume
(\ref{betalog}). Take $f\in \mathcal D^2_{0,0}$,
$f(z)=\sum_{n=1}^\infty a_nz^n$ ($z\in \mathbb D$). Using
(\ref{betalog}) and the fact that $\max (j,m)\,\asymp\,(j+m+1)$, we
see that
\begin{align*}\Vert \mathcal C_{(\eta )}(f)\Vert _{\mathcal
D^2_\beta }^2 \lesssim \,&\sum_{n=1}^\infty n^{1-\beta }\vert \eta
_n\vert ^2\left (\sum_{k=1}^n\vert a_k\vert \right )^2\\ =\,&
\sum_{n=1}^\infty n^{1-\beta }\vert \eta _n\vert ^2\left
(\sum_{j,m=1}^n\vert a_j\vert \vert a_m\vert \right )\\ \lesssim\,&
\sum_{j,m=1}^\infty \left (\vert a_j\vert \vert a_m\vert
\sum_{n=j+m+1}^\infty n^{1-\beta }\vert \eta _n\vert ^2\right )
\\ \lesssim\,&\sum_{j,m=1}^\infty \frac{\vert a_j\vert \vert
a_m\vert }{\log (j+m+1)}.
\end{align*}
Using Lemma\,\@2 of \cite{Shields} and Schur test as in
\cite[p.\,\@34]{GG-2025}, we deduce that $$\sum_{j,m=1}^\infty
\frac{\vert a_j\vert \vert a_m\vert }{\log (j+m+1)}\lesssim \Vert
f\Vert _{\mathcal D^2_{0,0}}^2$$ and, hence, $\mathcal C_{(\eta
)}\in \mathcal B(\mathcal D^2_{0,0}, \mathcal D^2_\beta )$.
\end{pf}
\par\bigskip

\centerline{\it Proof of  Theorem\,\@\ref{alphazero}\,\@(b).}
\par\medskip Assume that $\mathcal C_{(\eta )}\in \mathcal C(\mathcal
D^2_0, \mathcal D^2_\beta )$.  For $\frac{1}{2}<b<1$, let $h_b$ be
the function defined in the proof of Theorem\,\@\ref{D200beta}.
Clearly, $h_b\,\rightarrow\,0$, uniformly in compact subsets of
$\mathbb D$, as $b\to 1$, and $\Vert h_b\Vert _{\mathcal
D^2_0}^2\asymp 1$. Then
 Lemma\,\@3.\,\@7 of \cite{Tj} yields that $\Vert \mathcal C_{(\eta
 )}(h_b)\Vert _{\mathcal D^2_\beta }^2\to 0$, as $b\to 1$. We
 have seen above that for $b_N=1-\frac{1}{N}$ we have $$(1-\frac{1}{N})^{2N}\log N \sum_{n=N}^\infty
 n^{1-\beta }\vert \eta_n\vert ^2\,\lesssim \Vert \mathcal C_{(\eta
 )}(h_{b_N})\Vert _{\mathcal D^2_\beta }^2,\quad N\ge 2.$$
Therefore we obtain (\ref{opbetalog}) .
\par\medskip
Assume now that (\ref{opbetalog}) holds. Set
$$A_N\,=\,\sup_{j\ge N}(\log j)\sum_{n=j}^\infty n^{1-\beta }\vert
\eta_n\vert ^2,\quad N\ge 2.$$ Then $\{ A_N\} \to 0$, as $N\to
\infty $. Part (a) of the theorem gives that $\mathcal C_{(\eta
)}\in \mathcal B(\mathcal D^2_0, \mathcal D^2_\beta )$. As in the
proof of  Theorem \ref{alphaneg}, it suffices to prove that
$\mathcal C_N\to \mathcal C_{(\eta )}$ in the operator norm,
 where
$$\mathcal C_N(f)(z)=\sum_{n=0}^N\eta _n\left (\sum _{k=0}^na_k\right
)z^n,\quad z\in \mathbb D.$$ So, take $f\in \mathcal D^2_0$,
$f(z)=\sum_{n=0}^\infty a_nz^n$ ($z\in \mathbb D$). We have
$$\Vert \mathcal C_N(f)\,-\,\mathcal C_{(\eta )}(f)\Vert _{\mathcal
D^2_\beta }^2\,\le I\,+\,II\,+\,III,$$ where
\begin{align}\label{123}I=\vert \a_0\vert
^2\sum_{n=N+1}^\infty n^{1-\beta }\vert \eta _n\vert ^2,\quad &
II=\sum_{n=N+1}^\infty n^{1-\beta }\vert \eta_n\vert ^2\left
(\sum_{k=1}^N\vert a_k\vert\right )^2,\\ III=\sum_{n=N+1}^\infty
n^{1-\beta }&\vert \eta_n\vert ^2\left (\sum_{k=N+1}^n \vert
a_k\vert\right )^2.\nonumber\end{align} We have
\begin{align*}I=\vert a_0\vert ^2\sum_{n=N+1}^\infty n^{1-\beta }\vert
\eta_n\vert ^2\,\lesssim &\Vert f\Vert_{\mathcal D^2_0}^2(\log
N)\sum_{n=N}^\infty n^{1-\beta }\vert \eta_n\vert ^2\,\lesssim
A_N\Vert f\Vert_{\mathcal D^2_0}^2.
\end{align*}
\begin{align*}II\,=\,&\left(\sum_{n=N+1}^\infty n^{1-\beta }\vert \eta_n\vert
^2\right)\left (\sum_{k=1}^N\frac{k^{1/2}\vert a_k\vert }{k^{1/2}}\right)^2\\
\,\lesssim\,&\left(\sum_{n=N}^\infty n^{1-\beta }\vert \eta_n\vert
^2\right)\left (\sum_{k=1}^N\frac{1}{k}\right )\left
(\sum_{k=1}^\infty k\vert a_k\vert ^2\right )\,\lesssim A_N\Vert
f\Vert _{\mathcal D^2_0}^2.
\end{align*}
Arguing as in the proof of Theorem~\ref{D200beta}, we have
\begin{align*}III\,\le &\sum_{n=N+1}^\infty n^{1-\beta }\vert \eta
_n\vert ^2\left (\sum_{k=N+1}^n \vert a_k\vert \right )^2\\
\,\le &\sum_{j,m=N}^\infty \left (\vert a_j\vert \vert a_m\vert
\sum_{n=\max\{j,m\}}^\infty n^{1-\beta }\vert \eta _n\vert ^2\right
)
\\
 \,\lesssim
\, &\sum_{j,m=N}^\infty \left (\frac{\vert a_j\vert \vert a_m\vert
}{\log(\max\{j,m\})}A_{\max\{j,m\}}\right )
\\ \,\lesssim\,
&\,A_N\sum_{j,m=N}^\infty \left (\frac{\vert a_j\vert \vert a_m\vert
}{\log(j+m+1)}\right )\\ \,\lesssim\, & A_N\Vert f\Vert_{\mathcal
D^2_0}^2.
\end{align*} Putting everything together, we obtain that $\Vert
\mathcal C_N(f)-\mathcal C_{(\eta )}(f)\Vert _{\mathcal D^2_\beta
}^2\lesssim A_N\Vert f\Vert _{\mathcal D^2_0}$. Hence, $\mathcal
C_N\to \mathcal C_{(\eta )}$ in the operator norm.

\par\bigskip

\centerline{\it Proof of Theorem\,\@\ref{alphapos}.}\par\medskip Let
us start by proving the necessity of conditions (\ref{Oalphabeta})
and (\ref{oalphabeta}), for the boundedness and  the compactness,
respectively. For $\frac{1}{2}<b<1$, set
$$g_{b, \alpha }(z)\,=\,(1-b)^{\alpha /2}\sum_{n=1}^\infty n^{\alpha-1}b^nz^n,\quad z\in \mathbb D.$$ Then $\Vert g_{b, \alpha }\Vert _{\mathcal
D^2_\alpha }^2\,\asymp \,1$ and, since $\alpha >0$ it is clear that
$g_{b, \alpha }\to 0$, as $b\to 1$, uniformly in compact subsets of
$\mathbb D$.
\par For $N\ge 2$, we have
\begin{align*}\Vert \mathcal C_{(\eta )}(g_{b, \alpha })\Vert _{\mathcal D^2_\beta
}^2\,&\gtrsim (1-b)^\alpha \sum_{n=1}^\infty n^{1-\beta
}|\eta_n|^2\left
(\sum_{k=1}^nk^{\alpha -1}b^k\right )^2\\
&\gtrsim (1-b)^\alpha
\sum_{n=N}^\infty n^{1-\beta }|\eta_n|^2\left(\sum_{k=1}^N k^{\alpha-1}b^k\right)^2\\
&\gtrsim (1-b)^\alpha \sum_{n=N}^\infty n^{1-\beta }|\eta_n|^2N
^{2\alpha} b^{2N} &
\end{align*} Taking $b=1-\frac{1}{N}$,
we obtain
$$\Vert \mathcal C_{(\eta )}(g_{b, \alpha })\Vert _{\mathcal D^2_\beta
}^2\,\gtrsim N^{\alpha }\sum_{n=N}^\infty n^{1 -\beta }|\eta_n|^2.$$
\par If $\mathcal C_{(\eta )}\in \mathcal B(\mathcal D^2_\alpha ,
\mathcal D^2_\beta )$ then $\Vert \mathcal C_{(\eta )}(g_{b, \alpha
})\Vert _{\mathcal D^2_\beta }^2\lesssim 1$, and (\ref{Oalphabeta})
follows. If $\mathcal C_{(\eta )}\in \mathcal K(\mathcal D^2_\alpha
, \mathcal D^2_\beta )$ then $\Vert \mathcal C_{(\eta )}(g_{b,
\alpha })\Vert _{\mathcal D^2_\beta }^2\to 0$, as $b\to 1$,  and
(\ref{oalphabeta}) follows.
\par\medskip Let us turn to  prove the other implications.
We are going to use a result of G.~Bennett. Taking $r=2$ and $s=1$
in Theorem\,\@1 of \cite{Ben}, we obtain the following.
\begin{other}[Bennett\,\@\cite{Ben}]\label{ThBen} Let
$\{ u_n\} _{n=1}^\infty $, $\{ v_n\} _{n=1}^\infty $, and $\{ w_n\}
_{n=1}^\infty $ be three sequences of positive numbers. If
\begin{equation}\label{bigoh1}\sum_{n=1}^Nu_n\left (\sum_{k=1}^nv_k\right )^2\,=\,\og \left (
\,\sum_{n=1}^Nv_n\right ), \end{equation} then
\begin{equation}\label{bigoh2}\sum_{n=1}^Nu_n\left (\sum_{k=1}^nv_kw_k\right )^2\,=\,\og \left (
\sum_{n=1}^Nv_nw_n^2\right ).\end{equation}
\end{other}
\par\medskip
Assume that $\alpha >0$ and $\beta \in \mathbb R$ and that
(\ref{Oalphabeta}) holds.  In particular $F_{(\eta)}\in \mathcal
D^2_\beta$, that is $\sum_{n=1}^\infty n^{1-\beta }\vert \eta
_n\vert ^2\,<\,\infty $. \par\medskip It is well known (see
\cite[page 101]{Du:Hp}) that   (\ref{Oalphabeta}) is actually
equivalent to
 \begin{equation}\label{Oalphabeta1}
\sum_{n=1}^N n^{1+2\alpha -\beta }\vert \eta _n\vert^2\,=\,\og \left
(N^\alpha \right ),\quad \text{as $N\to \infty $.}\end{equation}

Take $f\in\mathcal D^2_{\alpha ,0}$, $f(z)=\sum_{n=1}^\infty a_nz^n$
($z\in \mathbb D$). We have
\begin{align}\label{uvw}\Vert \mathcal C_{(\eta )}(f)\Vert
_{\mathcal D^2_\beta }\,\lesssim\sum_{n=1}^\infty n^{1-\beta }\vert
\eta_n\vert ^2\left (\sum_{k=1}^n\vert a_k\vert \right
)^2\,=\,\sum_{n=1}^\infty u_n\left (\sum_{k=1}^nv_kw_k\right
)^2,\end{align} where,
$$u_k\,=\,k^{1-\beta }\vert \eta _k\vert ^2,\,\,\,v_k\,=\,k^{\alpha
-1}\,\,\,w^k\,=\,\frac{\vert a_k\vert }{k^{\alpha -1}},\quad k=1, 2,
3, \dots .$$ Using (\ref{Oalphabeta1}), we see that
\begin{align*}&\sum_{n=1}^Nu_n\left (\sum_{k=1}^nv_k\right
)^2\,=\,\sum_{n=1}^Nn^{1-\beta }\vert \eta _n\vert ^2\left
(\sum_{k=1}^nk^{\alpha -1}\right )^2\\
\,\lesssim\,&\sum_{n=1}^Nn^{1+2\alpha -\beta }\vert \eta _n\vert
^2\,\lesssim N^\alpha \,\asymp \,\sum_{n=1}^Nv_n.\end{align*} Then
Theorem\,\@\ref{ThBen} yields
\begin{align*}\sum_{n=1}^Nn^{1-\beta }\vert \eta _n\vert ^2\left
(\sum_{k=1}^n\vert a_k\vert\right )^2\,\lesssim
\sum_{n=1}^Nn^{\alpha -1}\frac{\vert a_n\vert ^2}{n^{2(\alpha -1)}}
\,=\,\sum_{n=1}^Nn^{1-\alpha }\vert a_n\vert ^2\,\lesssim \Vert
f\Vert_{\mathcal D^2_\alpha }^2.\end{align*} Thus it follows that
$\mathcal C_{(\eta )}\in \mathcal B\left (\mathcal D^2_{\alpha ,0},
\mathcal D^2_\beta \right )$. Since $\sum_{n=1}^\infty n^{1-\beta
}\vert \eta _n\vert ^2<\infty$, using Proposition\,\@\ref{X0} we
deduce that $\mathcal C_{(\eta )}\in \mathcal B\left (\mathcal
D^2_\alpha , \mathcal D^2_\beta \right )$.
\par Theorem\,\@\ref{ThBen} remains true if we
substitute \lq\lq big Oh\rq\rq \, by \lq\lq little oh\rq\rq \, in
(\ref{bigoh1}) and (\ref{bigoh2}). Using this we deduce that
(\ref{oalphabeta}) implies that $\mathcal C_{(\eta )}\in \mathcal
K(\mathcal D^2_{\alpha ,0}, \mathcal D^2_\beta )$ and, clearly, this
gives that $\mathcal C_{(\eta )}\in \mathcal K(\mathcal D^2_{\alpha
}, \mathcal D^2_\beta )$.

\par\bigskip

\par\bigskip
\section{Some final remarks and some further results}\label{final}
\begin{remark}\label{alphage1}
It is well known (see \cite[Page 101]{Du:Hp}) that if $\{ a_n\} $ is
a sequence non negative numbers and $0<\alpha_1<\alpha_2$, then
$$\sum_{n=0}^N n^{\alpha_2}a_n=O(N^{\alpha_1}) \Longleftrightarrow  \sum_{n=N}^\infty a_n=O(N^{\alpha_1-\alpha_2}).$$
Hence choosing $a_n=n^{1-\beta} |\eta_n|^2$ and either
$\alpha_1=\alpha$ and  $\alpha_2=2\alpha$  or  $\alpha_1=1$ and
$\alpha_2=1+\alpha$  for $\alpha>0$ we have
$$\sum_{n=1}^N n^{2\alpha }n^{1-\beta} |\eta_n|^2=O(N^{\alpha}) \Longleftrightarrow
\sum_{n=N}^\infty n^{1-\beta} |\eta_n|^2=O(N^{-\alpha})
\Longleftrightarrow  \sum_{n=1}^N n^{2+\alpha-\beta} |\eta_n|^2=O(N)
.$$

Similar results apply with \lq\lq little oh\rq\rq\, condition
instead of \lq\lq Big oh\rq\rq .
\end{remark}

\par Consequently, if $\alpha >1$, $\beta >-1$ and $\mu $ is a
complex Borel measure in $\mathbb D$ our results coincides with
those obtained by Lin and Xie in \cite{LX} for these values of
$\alpha $ and $\beta $.  But actually the result holds for any
$\alpha>0$ and any $\beta\in \mathbb R$.
\begin{corollary}\label{Dirichlet} If $\alpha>0 , \beta \in \R$ and $(\eta
)=\{ \eta _n\} $ is a sequence of complex numbers then:
$$\mathcal C_{(\eta )}\in \mathcal B(\mathcal D^2_\alpha , \mathcal D^2_\beta
)\,\,\,\,\Leftrightarrow\,\,\,\,\sum_{n=1}^N n^{\alpha +2-\beta
}\vert \eta_n\vert ^2\,=\,\og (N).$$
$$\mathcal C_{(\eta )}\in \mathcal K(\mathcal D^2_\alpha , \mathcal D^2_\beta
)\,\,\,\,\Leftrightarrow\,\,\,\,\sum_{n=1}^N n^{\alpha +2-\beta
}\vert \eta_n\vert ^2\,=\,\op (N).$$
\end{corollary}
Also, Corollaries 3.\,\@6 and 3.\,\@7 of \cite{LX} can be extended.
Using that $A^2_\alpha= \mathcal D^2_{\alpha+2}$ for $\alpha>-1$ in
Theorem \ref{alphapos} and the previous remark we can write the
following corollary.
\begin{corollary}\label{Bergman} If $\alpha , \beta >-1$ and $(\eta
)=\{ \eta _n\} $ is a sequence of complex numbers then:
$$\mathcal C_{(\eta )}\in \mathcal B(A^2_\alpha , A^2_\beta
)\,\,\,\,\Leftrightarrow\,\,\,\,\sum_{n=1}^N n^{\alpha +2-\beta
}\vert \eta_n\vert ^2\,=\,\og (N) \,\,\,\,\Leftrightarrow\,\,\,\
\sum_{n=N}^\infty n^{-(\beta+1)}|\eta_n|^2= \og(N^{-(\alpha+2)}).$$
$$\mathcal C_{(\eta )}\in \mathcal K(A^2_\alpha , A^2_\beta
)\,\,\,\,\Leftrightarrow\,\,\,\,\sum_{n=1}^N n^{\alpha +2-\beta
}\vert \eta_n\vert ^2\,=\,\op (N) \,\,\,\,\Leftrightarrow\,\,\,\
\sum_{n=N}^\infty n^{-(\beta+1)}|\eta_n|^2= \op(N^{-(\alpha+2)}).$$
\end{corollary}
\par\bigskip

\begin{remark}\label{decreasing} If $\alpha >0$ and the sequence $\{ \vert
\eta_n\vert \} $ is decreasing then (\ref{Oalphabeta}) is equivalent
to
\begin{equation}\label{etaog} \vert \eta _N\vert \,=\,\og \left
(\frac{1}{N^{1+\frac{\alpha -\beta }{2}}}\right).\end{equation}
\end{remark}
Indeed, since $\{ \vert \eta_n\vert \} $ is decreasing, we see that
$$\vert \eta_{2N}\vert ^2N^{2 -\beta }\,\leq C \,\sum_{n=N}^{2N}n^{1 -\beta }\vert \eta_{n}\vert
^2\,\le \, C\sum_{n=N}^\infty n^{1- \beta }\vert \eta _n\vert ^2.$$
Then we see that (\ref{Oalphabeta}) implies (\ref{etaog}). The other
implication is obvious since $\alpha>0$.
\par\bigskip
\begin{remark}\label{decreasing2} If  $\alpha >0$ and the sequence $\{ \vert
\eta_n\vert \} $ is decreasing as in remark\,\@\ref{decreasing} and
$2+\alpha \le \beta $, then (\ref{etaog}) simply says that $\vert
\eta_N\vert \,=\, \og \left (1\right )$.
\end{remark}

\par\medskip  If $\mu $ is a finite and positive Borel measure on
$[0,1)$ and, for $n=0, 1, 2, \dots $, $\mu_n$ is the moment of order
$n$ of $\mu$, then the sequence $\{ \mu_n\} $ is a decreasing
sequence of non-negative numbers. Then Theorem\,\@\ref{alphapos},
Remark\,\@\ref{decreasing}, and Remark\,\@\ref{decreasing2} give the
following result that extends those in \cite{JT}  and
\cite{Blasco1}.
\begin{theorem}\label{alphapos-borel-measure} Suppose that $\alpha >0$
and let $\mu $ be a finite and positive Borel measure on $[0,1)$.
\begin{itemize}
\item[(a)] If $2+\alpha \le \beta $ then $\mathcal C_\mu \in \mathcal B(\mathcal D^2_\alpha , \mathcal
D^2_\beta )$. \item[(b)] If $\beta <\alpha +2$, then $\mathcal C_\mu
\in \mathcal B(\mathcal D^2_\alpha , \mathcal D^2_\beta )$ if and
only if (\ref{etaog}) holds.\end{itemize}\end{theorem}
\par\medskip
As a consequence of Theorem\,\@\ref{alphapos-borel-measure}\,\@(b),
we obtain.
\begin{theorem}\label{radialmeasure} If $\alpha >0$, $\beta <\alpha +2$,  and $\mu $ is a finite and positive Borel measure on
$[0,1)$, then $\mathcal C_\mu \in \mathcal B(\mathcal D^2_\alpha ,
\mathcal D^2_\beta )$ if and only it $\mu $ is a $1+\frac{\alpha
-\beta }{2}$-Carleson measure, that is,
$$\mu ([t, 1))\,\lesssim \,(1-t)^{1+\frac{\alpha
-\beta }{2}}.$$
\end{theorem}

\par The result follows using
Theorem\,\@\ref{alphapos-borel-measure}\,\@(b), the fact that
$1+\frac{\alpha -\beta }{2}>0$, and Proposition\,\@1 of
\cite{Ch-Gi-Pe}.
\par\bigskip
{\bf Data Availability.} All data generated or analyzed during this
study are included in this article and in its bibliography

\par\medskip
{\bf Conflict of interest.} The authors declare that there is no
conflict of interest.

\par
\end{document}